\documentclass[a4paper,10pt]{article}
\usepackage{latexsym}
\usepackage{amsmath}
\usepackage{amsfonts}
\usepackage{amssymb}
\usepackage{vmargin}
\usepackage{graphicx}
\usepackage{url}

\setmarginsrb{1.5cm}{1.5cm}{1.5cm}{2cm}{0mm}{0mm}{0mm}{10mm}
\usepackage{bbding}
\usepackage{color}

\usepackage[leftmargin=3em]{quoting}
\newcounter{examplecounter}
\newenvironment{example}
{\begin{quoting}[leftmargin=0pt]
\small
\refstepcounter{examplecounter}%
\textbf{Example \arabic{examplecounter}}%
\quad
}{%
\end{quoting}%
}

\usepackage{caption}
\usepackage[labelformat=simple]{subcaption}
\usepackage{calc}
\usepackage{enumitem}
\usepackage{kantlipsum}
\usepackage{array}
\usepackage{breakurl}

\newcommand{\Mark}[1]{\textsuperscript{#1}}

\usepackage{titlesec}

\begin{document}
\twocolumn[{%
 \centering
 \Large {\bf Embedded librarianship and problem-based learning in undergraduate mathematics courses} \\[1.5em]
 \normalsize N. Karjanto\Mark{1},
        M. Kairatbekkyzy\Mark{2},
        and J. Agee\Mark{3}\\[1em]
 \small
 \begin{tabular}{*{2}{>{\centering}p{.5\textwidth}}}
  \Mark{1}Department of Mathematics & \Mark{2}Nazarbayev University Library \tabularnewline
  School of Science and Technology  & Astana, Kazakhstan  \tabularnewline
  Nazarbayev University             & \Mark{3}American University of Phnom Penh Library \tabularnewline 
  \url{natanael.karjanto@gmail.com} & Phnom Penh, Cambodia  
 \end{tabular}\\[3em] 
}]

\begin{abstract}
{\footnotesize
A pedagogical approach of problem-based learning with embedded librarianship in several undergraduate mathematics courses is implemented in this educational research.
The students are assigned to work on several projects on various applications of mathematical topics in daily life and submit written reports.
An embedded librarian collaborates together with the instructor and the students to improve the students' information literacy.
Initial reaction and anecdotal evidence show that the students' information literacy and academic performance have improved throughout the semesters.}\vspace*{0.2cm}\\
\end{abstract}
{\setlength{\parindent}{0pt}
{\footnotesize {\slshape Keywords:} problem-based learning, undergraduate mathematics, embedded librarian, library, instructional technology, Kazakhstan \par}
}

\section{Introduction}

During an annual state of the national address in Kazakhstan on December 15, 2012,
the current president of Kazakhstan, Nursultan Nazarbayev, announced the `Kazakhstan 2050 Strategy'~\cite{2050}.
In general, this strategy drives for political, economic and social reforms to position
the country in the top 30 nations by the year 2050 in term of economic strength.
In addition, in the context of higher education, based on the president's initiative himself,
an autonomous research-based university has been founded in 2010~\cite{nu},
where the authors of this article are affiliated and this educational research has been carried out for the past one year.
This university and the related Nazarbayev Intellectual Schools are dedicated to promote educational reform in the country.
The current university's president, Dr. Shigeo Katsu, once indicated that these institutions are created outside the system in order to reform the current, existing system~\cite{katsu13}.

One aspect of the Kazakhstan 2050 Strategy that has a particular interest in the context of this article
is in the realm of (higher) education, in particular the modernization of teaching methods~\cite{2050}.
Motivated by this, as well as the fact that the traditional style of lecture-based teaching is far from effective~\cite{hanford12},
we are interested to shift the paradigm of higher education from traditional teacher-centered pedagogy to a student-centered active learning style.
One prime example of this non-traditional teaching style is to implement problem-based learning (PBL) coupled with embedded librarianship.
The focus of this educational research is on the implementation of PBL in several undergraduate mathematics courses served by
Department of Mathematics, School of Science and Technology (SST) at Nazarbayev University to the students at both SST and SHSS (School of Humanities and Social Sciences).
This type of educational research on this collaboration between an embedded librarian and an instructor for a good practice of PBL implementation
seems novel in the field of mathematics and science courses.
There are, however, some reports of embedded librarianship in action in other subjects too.
Some examples are in a freshman speech class~\cite{hall08}, in an English Composition class~\cite{jacobs09},
in health sciences at McMaster University in Canada~\cite{yates11} and the Ross School of Business at the University of Michigan~\cite{rupp11}.

This article is organized as follows.
The following section deals with embedded librarianship. Section~\ref{PBL} discusses briefly what is meant by the pedagogical approach of problem-based learning.
Some real classroom examples in several undergraduate mathematics courses are also presented.
Section~\ref{anecdote} follows with initial reactions from the students on the implementation of PBL with embedded librarianship.
Some anecdotal evidences, both qualitative and quantitative, are also presented in this section.
Finally, Section~\ref{conclusion} provides conclusion and remark to our observation.

\section{Embedded librarianship}

The onset of new technologies and the students' approaches to gain information have dramatically changed work methods in the academic world as well as in academic libraries.
If in the past the library's main role is to store the information generated by society and to provide people with access to this information with the help of a librarian,
the situation has changed recently.
With the invention of the Internet, the roles of libraries and librarians have also shifted from a physical setting to a more virtual environment.
The World Wide Web has assumed the library's role as the storage of information and a well-known search engine {\slshape Google} has assumed the librarians' role as an access provider to this information.
At this point, the phrase ``Do not wait until you are asked" is becoming relevant for academic librarians.
The librarians' role now is to focus on the users rather than the collection of information.
This includes instruction sessions, library services for teaching and learning, and support for teaching and research staff.

Embedded academic librarians find ways to embed their skills and services in physical and virtual environments~\cite{brower11}.
The author also describes six characteristics of embedded librarianships.
{\sl First}, they collaborate with their users.
{\sl Second}, they form partnerships on the department and campus levels.
{\sl Third}, they provide needs based services.
{\sl Fourth}, they offer convenient and user-friendly services outside of library settings.
{\sl Fifth}, they become immersed in the culture and spaces of users.
{\sl Sixth}, they understand the discipline including the culture and research habits of their users.


Furthermore, as stated in the Special Libraries Association report,
embedded librarianship involves focusing on the needs of one or more specific groups,
building relationships with these groups, developing a deep understanding of their work
and providing information services that are highly customized and targeted to their greatest needs.
In effect, it involves shifting the basis of library services from the traditional, transactional,
question-and-answer model of reference services to one in which there is high trust, close collaboration
and shared responsibility for outcomes~\cite{shumaker09}.

Motivated by the transformation in information services mentioned in the opening paragraphs of this section,
the Nazarbayev University Library has started its pilot project in Autumn/Fall 2013 to offer faculty the embedded librarian service.
This pilot project is conducted to determine whether the service is useful and practical for both faculty and students to support teaching and learning activities.
If it turns out positive, this collaboration service will be continuously implemented in the future and possibly to be expanded with more faculty being involved.

Although the practice of PBL have been conducted in several undergraduate mathematics courses, however,
PBL with strong embedded librarian involvement has been conducted particularly in two different courses,
i.e. Discrete Mathematics, and Linear Algebra with Applications, during Autumn/Fall 2013 and Spring 2014, respectively.
The embedded librarian plays a role in assisting the students to complete their PBL assignments and to improve their information literacy.
There are six steps implemented by the embedded librarian to assist the students in completing their PBL activities (also known as `embedded librarian in action').
{\bf Step~1:} Embedded librarian conducts two library in-class sessions.
The first session is the instruction on how to use library resources for a particular topic and
the second one is an instruction on how to work with {\slshape RefWorks} citation tool.
{\bf Step~2:} Providing embedded support in {\slshape Moodle} virtual learning environment.
This step includes providing information, instructional presentations on relevant resources and participation in an online chat.
{\bf Step~3:} Conducting individual tutorials, usually in the library.
This is provided for the students who have follow-up questions after the library sessions in the classroom,
and for those who might need more detailed explanation on how to progress in their PBL assignments.
{\bf Step~4:} Evaluation of reference list sources as soon as the PBL assignments are submitted by students.
{\bf Step~5:} Analysis of the students' progress in the use of scholarly resources and their ability to cite them correctly.
{\bf Step~6:} If necessary, another extra in-class session will be conducted with detailed explanations of
potential mistakes made by the students and discussion for better presentations.


\section{PBL in mathematics courses} \label{PBL}

Problem-based learning (PBL) is a student-centered pedagogy in which students learn about a subject through the experience of problem solving.
It was pioneered by Howard Barrows and his colleagues in the late 1960s at the medical school program at McMaster University in Hamilton, Ontario, Canada.
The features of PBL include the students to enhance their ability to learn both thinking strategies and domain knowledge, active learning style, small group collaborative learning,
the role of instructor as facilitator: from a sage on a stage to a guide on the side~\cite{king93}, a problem forms a basis for learning and the problem stimulates cognitive process~\cite{barrows96}.
The are a number of evidences where PBL model is more superior than the traditional teaching style.
For example, Slovenian students exposed to PBL are better at solving more difficult mathematics problems than their counterparts~\cite{cotic09}.
Middle-school's performance in standardized tests is also improving after PBL is implemented~\cite{hmelo07}.
Physicians' competence after graduation, including communications skills, are retained~\cite{koh08}.

However, there are also criticism toward the PBL style. Apart from time consuming~\cite{hung11}, there is a heavy cognitive load for the learners~\cite{sweller88,sweller06}.
It is also found that PBL is less effective instructional strategy than studying worked examples~\cite{sweller85,cooper87},
where the authors conducted several classroom studies where students studying algebra problems.
Implementing PBL requires lots of resources, planning and organization~\cite{azer11a,azer11b}.
The author also mentioned 12 steps that need to be implemented for a successful PBL practice.
It is shown that there is no convincing evidence that PBL improves knowledge base and clinical performance, at least not of
the magnitude that would be expected given the resources required for a PBL curriculum~\cite{colliver00}.

Some PBL examples are given as follows.
\begin{example} \label{ex:astanaclean} (Keep Astana Clean!)\\
You are working as a consultant team for the municipal services of Astana city.
During the summer time, sometimes strong wind carries dust storms into the city, while during the winter time, the snow envelopes the city.
The city council wants to keep our beautiful city Astana to be clean throughout the year, in particular from the dust (during the summer) and the snow (during the winter).
Since the city councils would like to save millions of tenge annually to clean the city, they consult your team on minimizing the cost.
You and your team do not want to disappoint your boss, so your team proposes some strategies and reports them to municipal services of Astana regarding the best approach to handle the situation.
\end{example}

\begin{example} \label{ex:lotka} (Predator-prey system) \\
Explore one of the following relationships (choose only one) between populations of predators and preys:
sharks and food fish, ladybugs and aphids, wolves and rabbits, and other animals that you find interesting.
The related material is a chapter on Introduction to Differential Equations. Model the interaction with differential equations.
You may approach the problem from analytical, numerical, and graphical viewpoints.
Conclude your findings. Provide a well-documented bibliography.
The bibliography should contain a minimum of 10 peer-reviewed journal articles.
Follow APA (American Psychological Society) style for reference list and in-text citation. See \url{www.apastyle.org}.
A good report should contain and explain, but not limited to, the followings:
mathematical model, equilibrium points, conserved quantity and conservation law, lags behind, some plots.
\end{example}



\begin{example} \label{ex:leontief} (Leontief's model for economy)\\ 
Economy of a country is a complicated network of interdependence. Changes in one part can affect other parts.
Economists (and mathematicians) have attempted to describe and to make predictions of such complicated relationship using a mathematical model of system of linear equations.
Based on Leontief's input-output model, describe the economy of Kazakhstan using the linear system that you have learned in the class.
Make a maximum five-page report as a three-member team. This problem is ill-structured, you need to assume something.
Provide a well-documented bibliography to support your report. Consult the subject librarian if necessary.
\end{example}

\section{Some results} \label{anecdote}

\subsection{Initial reaction}
Implementing good PBL practices in classroom environments is very challenging.
Apart from time consuming, both inside and outside the classroom settings,
they require plenty of preparation from the instructors' side.
In addition, the students do not usually expect to solve non-routine problems during class sessions.
When the students attending classes, generally they prefer to sit down and listen to the instructor who would tell them something.
Furthermore, when the instructor writes essential information on the board, the students take notes.
This is typically how the traditional-conventional classroom looks. The instructor often plays as a sage on the stage.

When the PBL assignments are distributed, some students show various reactions.
The students may be confused about where and how to start the project. Some of them are perplexed and doubt whether they have made particular assumptions in the right direction.
They do not expect to conduct a project assignment in a mathematics class, which is usually geared toward problem-solving skills rather than
conducting a literature study and writing a research report.
In particular, the students majoring in Economics with mathematics emphasis have expressed their concern that
they have already many writing assignments from other humanities and social science classes, so they prefer to have only practice sessions on problem-solving instead of PBL sessions.

In connection to working cooperatively with the embedded librarian, the initial reaction from the students is quite positive.
The students through their team representatives have approached the embedded librarian to utilize library resources, to locate and select information,
and to include relevant references in their project reports.

\subsection{Anecdotal evidence}

\subsubsection{Qualitative comparison}

A qualitative comparison between Discrete Mathematics and Linear Algebra with Applications classes due to embedded librarianship is presented in Table~\ref{qualitative}.
This comparison is conducted between two different groups taking two different classes in two different semesters where embedded librarian is extensively involved in grading the submitted PBL assignments. The classes are Discrete Mathematics and Linear Algebra with Applications, conducted in Autumn/Fall 2013 and Spring 2014, respectively.
It is observed that there was significant improvement between the former group and the latter one in terms of reference list contents and citation style.
Some possible explanations for the improvement are clearer specific instructions and stronger imposed remarks on these aspects.
\begin{table}[t]
\begin{center}
{\footnotesize
\begin{tabular}{|l|l|l|}
\hline \hline
Course     & Discrete Mathematics & Linear Algebra \\ \hline
Semester   & \quad Autumn/Fall 2013   & \qquad Spring 2014 \\ \hline
Students   & \qquad \qquad 18     & \qquad \qquad  29 \\ \hline
Team       & \qquad \qquad \;6   & \qquad \qquad \;9 \\ \hline
           & An average number    & An average number \\
Reference  & of references $= 5$ & of references $= 10$\\
list       & Source: textbooks   & Library e-resources  \\
contents   & internet resources     & Government docs  \\
           & No library e-resources  & Reliable statistics \\ \hline
Citation   & No style               & APA style  \\
style      & No proper citation & Few proper citation  \\
           & Difficult to follow   & Easier to follow  \\ \hline
Evaluation & Excellent or good: & Excellent or good: \\
           & {\scriptsize \color{blue} \FiveStar}{\scriptsize\color{red}\FiveStarOpen \FiveStarOpen \FiveStarOpen \FiveStarOpen \FiveStarOpen \FiveStarOpen}
           & {\scriptsize \color{blue} \FiveStar \FiveStar \FiveStar \FiveStar \FiveStar}{\scriptsize\color{red}\FiveStarOpen \FiveStarOpen \FiveStarOpen \FiveStarOpen} \\  \hline \hline
\end{tabular}
}
\caption{\footnotesize Qualitative comparison between Discrete Mathematics and Linear Algebra classes due to embedded librarianship.}
\label{qualitative}
\vspace*{-0.75cm}
\end{center}
\end{table}
\begin{table}[b]
\vspace*{-0.25cm}
\begin{center}
{\footnotesize
\begin{tabular}{|c|c|c|c|c|c|c|}
\hline \hline
        & \multicolumn{3}{|c|}{Calculus~2} & \multicolumn{3}{|c|}{Linear Algebra} \\ \hline
Group   & N  & Mean  & SD    & N  & Mean  & SD    \\ \hline
Non-PBL & 34 & 73.72 & 10.37 & 29 & 63.67 & 18.17 \\ \hline
PBL     & 30 & 84.42 & 10.32 & 29 & 75.85 & 14.21 \\ \hline \hline
                      & \multicolumn{3}{|c|}{Calculus~2}    & \multicolumn{3}{|c|}{Linear Algebra} \\ \hline
DF    & \multicolumn{3}{|c|}{61}          & \multicolumn{3}{|c|}{53}    \\ \hline
$|t|$  & \multicolumn{3}{|c|}{4.130}       & \multicolumn{3}{|c|}{2.861} \\ \hline
$p$-value             & \multicolumn{3}{|c|}{0.0001}      & \multicolumn{3}{|c|}{0.006} \\ \hline
Conclusion            & \multicolumn{3}{|c|}{very highly significant} & \multicolumn{3}{|c|}{highly significant} \\
\hline \hline
\end{tabular}
}
\caption{\footnotesize Quantitative comparison between non-PBL and PBL groups in two mathematics courses: Calculus~2 and Linear Algebra.
DF indicates the degrees of freedom and $|t|$ means the statistic value correspond to two-sample Welsch's $t$ test.}
\label{quantitative}
\end{center}
\end{table}

\subsubsection{Quantitative comparison}
A quantitative comparison between PBL and non-PBL groups in two sections of two mathematics courses Calculus~2 and Linear Algebra with Applications is presented in Table~\ref{quantitative}.
For two Calculus~2 sections, the PBL and non-PBL groups are both enrolled in Spring 2014.
For Linear Algebra classes, the non-PBL and PBL groups are enrolled in Spring 2013 and Spring 2014, respectively.
The statistical analysis is conducted using the two-sample Welsch's $t$ test with an assumption that the populations differ in variance ($\sigma_1^2 \neq \sigma_2^2$).
The means and the standard deviations for both classes are taken from the averages of the two midterm tests and a final exam, where the weight for the final exam accounts twice of the midterm tests.
It is obtained that the comparisons between two groups in two different classes of Calculus~2 and Linear Algebra with Applications, the $p$-values are very small, 0.0001 and 0.006, respectively.
These values suggest that the difference between non-PBL and PBL groups in terms of academic performance, in this context indicated by the exam results, are very highly significant and highly significant, respectively, for which the PBL groups outperformed the non-PBL ones in both classes.

\section{Conclusion} \label{conclusion}
An educational research on shifting of higher education paradigm from a traditional-didactic teacher-centered teaching style
to a non-conventional student-centered pedagogy in several undergraduate mathematics courses has been presented in this article.
The approach is by implementing problem-based learning (PBL) activities with embedded librarianship, where the instructor, the embedded librarian, and the students work collaboratively
in completing project assignments related to topics covered in particular classes.
The PBL examples are designed with local contexts and applications in mind.
They are challenging in nature and most of them are ill-structured.
Thus, the students need to brainstorm, discuss, arrive to particular assumptions, and find relevant literature in order to write a good PBL report.
Even though some initial reactions indicate that this approach is quite demanding and challenging for all parties,
some anecdotal evidences show that the students' information literacy, writing skills and academic performances have improved significantly.
More extended qualitative and quantitative comparisons need to be conducted for future research to confirm the effectiveness of PBL practice with embedded librarianship.

\section*{\large Acknowledgement}
{\footnotesize The authors would like to acknowledge Professor Carole Faucher (Graduate School of Education, Nazarbayev University, Kazakhstan),
Dr. Loretta O'Donnell (Vice Provost Academic Affairs, Nazarbayev University, Kazakhstan),
Dr. Rob Lahaye (Department of Physics and University College, Sungkyunkwan University, South Korea),
Dr. Tu\u{g}rul Kar (Department of Primary Mathematics Education, Atat\"{u}rk University, Turkey),
Elizabeth Yates (James A. Gibson Library, Brock University, Canada) and
Nathan Rupp (Cushing/ Whitney Medical Library, Yale University, Connecticut, USA)
for many fruitful discussions and constructive suggestions for the improvement of this paper as well as
Professor Sang-Gu Lee (Department of Mathematics, Sungkyunkwan University, South Korea) as the
scientific and local organizing committee of the 19th International Linear Algebra Society (ILAS 2014) conference,
Dr. Ajit Kumar (Department of Mathematics, Institute of Chemical Technology, Mumbai, India) and
Professor Abraham Berman (Department of Mathematics, Technion--Israel Institute of Technology, Israel) as the organizers of Contributed Mini-Symposia (CMS-4) of ILAS 2014,
Professor Yong-Ju Choie (Department of Mathematics, POSTECH, South Korea) and
Professor Sun-Sook Noh (Department of Mathematics Education, Ewha Womans University, South Korea) as the chairpersons of the local program and local organizing committee, respectively,
of International Congress of Women Mathematicians ICWM 2014.
\par}


{\footnotesize
\begin{thebibliography}{99}
\bibitem{azer11a} Azer, S. A. (2011) 
\textit{Medical Education} {\bf 45}: 510.

\bibitem{azer11b} Azer, S. A. (2011) 
\textit{Medical Teacher} {\bf 33}(10): 808–-813.

\bibitem{barrows96} Barrows, H. S. (1996) 
\textit{New Directions for Teaching and Learning} {\bf 68}: 3--12.

\bibitem{brower11} {Brower, M. (2011) 
In Kvenild, C. and Calkins, K. (Eds.) \textit{Embedded Librarians: Moving beyond one-shot instruction}. Chicago: Association of College and Research Libraries. \par}

\bibitem{christensen08} Christensen, O. L. (2008) 
\textit{Teaching Mathematics and Its Applications} {\bf 27}: 131--139.

\bibitem{colliver00} Colliver J. A. (2000) 
\textit{Academic Medicine}. {\bf 75}: 259--266.

\bibitem{cooper87} Cooper, G. A. and Sweller, J. (1987) 
\textit{Journal of Educational Psychology} {\bf 79}: 347--362.

\bibitem{cotic09} Coti\v{c}, M. and Zuljan, M. V. (2009) 
\textit{Educational Studies} {\bf 35}: 297--310.

\bibitem{crosby00} Crosby, O. (2000) 
\textit{Occupational Outlook Quarterly} {\bf 44}: 3--15.

\bibitem{hall08} Hall, R. A. (2008) 
\textit{C\&RL News} {\bf 69}(1): 28--30.

\bibitem{hanford12} Hanford, E. (2012) 
\textit{National Public Radio NPR}: 4~pp. Available online at \url{http://www.npr.org/2012/01/01/144550920/physicists-seek-to-lose-the-lecture-as-teaching-tool}. Last accessed on June 11, 2014.

\bibitem{hmelo07} Hmelo-Silver, C. E., Duncan, R. G. and Chinn, C. A. (2007) 
\textit{Educational Psychologist} {\bf 42}: 99-107.

\bibitem{hung11} Hung, W. (2011) 
\textit{Educational Technology Research and Development} {\bf 59}(4): 529--552.

\bibitem{isik11} I\c{s}ik, A., Kaplan, A., I\c{s}ik, C., Konyalio\u{g}lu, A. C., G\"{u}ler, G. and Kar, T. (2011) 
\textit{International Journal of Humanities and Social Science} {\bf 1}: 187--196.

\bibitem{jacobs09} Jacobs, H. L. M. and Jacobs, D. (2009) 
\textit{References \& User Services Quarterly} {\bf 49}(1): 72--82.

\bibitem{katsu13} Katsu, S. (2013) Speech delivered during onboarding meeting for Nazarbayev University employees.


\bibitem{2050} Kazakhstan 2050 Strategy's official website (2014). Available online at \url{http://strategy2050.kz/en/}.

\bibitem{king93} King, A. (1993). 
\textit{College Teaching} {\bf 41}(1): 30--35.

\bibitem{koh08} Koh, G. C.-H., Khoo, H. E., Wong, M. L. and Koh, D. (2008) 
\textit{Canadian Medical Association Journal} {\bf 178}: 34--41.

\bibitem{nu} Nazarbayev University's official website (2014). Available online at \url{http://nu.edu.kz}.

\bibitem{regalado10a} Regalado-M\'{e}ndez, A., B\'{a}ez-Gonz\'{a}lez, J. G., Peralta-Reyes, E. and Trujillo-Tapia,  M. N. (2010)
\textit{Technological Developments in Networking, Education and Automation}: 13--18.

\bibitem{regalado10b} Regalado-M\'{e}ndez, A., Cid-Rogr\'{i}guez, M. R. and B\'{a}ez-Gonz\'{a}lez, J. G. (2010)
\textit{International Journal of Software Engineering and Applications} {\bf 1}: 54--73.

\bibitem{rupp11} Rupp, N. (2011) 
\textit{Stephen M. Ross School of Business, Kresge Business Library--Papers \& Presentation Series}: 6~pp. Available online at \url{http://deepblue.lib.umich.edu/bitstream/handle/2027.42/86726/Rupp_Action_Based_Learning_SLA_Michigan_Newsletter.pdf?sequence=1}. Last accessed on June 11, 2014.

\bibitem{shumaker09} Shumaker, D. and Talley, M. (2009) \textit{Models of embedded librarianship: Final report}, {\bf 9}: 1--195. Special Libraries Association.
Available online at \url{http://hq.sla.org/pdfs/EmbeddedLibrarianshipFinalRptRev.pdf}. Last accessed on June 6, 2014.

\bibitem{sweller88} Sweller, J (1988) 
\textit{Cognitive Science} {\bf 12}(2): 257--285.

\bibitem{sweller06} Sweller, J. (2006) 
\textit{Learning and Instruction} {\bf 16}: 165--169.

\bibitem{sweller85} Sweller, J. and Cooper, G. A. (1985) 
\textit{Cognition and Instruction} {\bf 2}: 59--89.

\bibitem{yates11} Yates, E. (2011) 
University of Western Ontario, Guided research 9411: 52~pp. Available online at \url{http://fims-grc.fims.uwo.ca/wp-content/uploads/2013/06/LIS9411_Yates.pdf}. Last accessed on June 11, 2014.
\end{thebibliography}
}
\end{document}